\documentclass[journal]{IEEEtran}

\usepackage[utf8]{inputenc}
\usepackage{bbm}

\usepackage[final]{hyperref} % adds hyper links inside the generated pdf file

\newcommand{\E}{\mathbb{E}}
\newcommand{\Lop}{\mathrm{L}}
\newcommand{\Dop}{\mathrm{D}}

\newcommand{\Sp}{\mathcal{S'}(\mathbb{R})}
\newcommand{\Swa}{\mathcal{S}(\mathbb{R})}
\usepackage{amsthm}
\usepackage{amsfonts}
\usepackage{color}
\usepackage{mathtools, stmaryrd}
\usepackage{tikz}
\usetikzlibrary{automata, positioning}
\usepackage{graphicx}
\usepackage[]{algorithm2e}
\usepackage{svg}
\usepackage{amssymb}
\usepackage{multirow}
\newtheorem{mydef}{Definition}

\newtheorem{prop}{Proposition}
\usepackage[]{algorithm2e}
\usepackage{cite}
\usepackage{url}

\begin{document}

\title{Generating Sparse Stochastic Processes   Using Matched Splines}

\author{Leello~Dadi\thanks{Leello Dadi and Shayan Aziznejad contributed equally to this work.},
        Shayan Aziznejad, 
        and~Michael~Unser,~\IEEEmembership{Fellow,~IEEE}
\thanks{This work was done at the Biomedical Imaging Group, \'{E}cole polytechnique f\'{e}d\'{e}rale de Lausanne,  Switzerland (e-mail: leello.tadesse@gmail.com; shayan.aziznejad@epfl.ch; michael.unser@epfl.ch). It was funded by the Swiss National Science Foundation under Grant 200020\_184646 / 1.}
 }

% The paper headers
%\markboth{Journal of \LaTeX\ Class Files,~Vol.~14, No.~8, August~2015}%
%{Shell \MakeLowercase{\textit{et al.}}: Generating sparse stochastic %processes}

% make the title area
\maketitle
\thispagestyle{plain}
\pagestyle{plain}

\begin{abstract}
We provide an algorithm to generate trajectories of sparse stochastic processes that are solutions of linear ordinary differential equations driven by Lévy white noises. A recent paper showed that these processes are limits in law of generalized compound-Poisson processes. Based on this result, we derive an off-the-grid algorithm that generates  arbitrarily close approximations of the target process. Our method relies on a \mbox{B-spline} representation of generalized compound-Poisson processes. We illustrate numerically the validity of our approach. 
\end{abstract}


\IEEEpeerreviewmaketitle

\section{Introduction}

Motivated by tractability and results such as the central-limit theorem, most of the early work in statistical signal processing has focused on Gaussian models \cite{gray2004introduction}. In particular, the theory of Gaussian stationary processes provided justifications for the use of the discrete cosine transform \cite{ahmed_discrete_1974}  as an approximation of the Karhunen–Loève transform, and  the Kalman filter \cite{kalman_new_1960} as an optimal estimator.

However, the analysis of real-world signals   has revealed that the Gaussian framework may be insufficient to capture the breadth of the underlying behaviors \cite{mumford_pattern_2010,srivastava2003advances}.  An important property that escapes the Gaussian framework is that of sparsity in some transform domains\cite{amini_compressibility_2011}. Sparsity being an essential component of modern signal processing \cite{mallat_wavelet_1999, donoho_compressed_2006, starck_sparse_2010}, the authors of \cite{unser_introduction_2014}   proposed a wider stochastic framework that encompasses both Gaussian and sparsity-compatible models. Within this framework, a continuous-time signal is a realization of a stochastic process $s$ that can be whitened by some linear, shift-invariant operator $\mathrm{L}$. The key here is that the resulting white noise, or {innovation}, is not necessarily Gaussian. Put formally, signals are solutions of

\begin{equation}
    \mathrm{L}s = w,
\label{eq::found}
\end{equation}
where $w$ is a well defined innovation process called a {Lévy white noise} \cite{kailath1970innovations}. The term \textit{Lévy} here comes from the fact that $w$ is an object that can be interpreted as the derivative of a Lévy process in the sense of distributions \cite{sato, unser_unified_2014}. Whenever $w$ is non-Gaussian, the realizations of $s$ can be shown to be sparse. Accordingly, they have been named sparse stochastic processes  \cite{unser_introduction_2014}. Specific instances of such processes have been  used to model  natural signals  such as images \cite{mumford2001stochastic,bostan_sparse_2013}, RF echoes in ultrasound \cite{kutay2001modeling}, and network traffic in communication systems \cite{kogon1996signal,gallardo2000use,laskin2002fractional}.
% and, DNA sequences in computational biology \cite{allegrini1995dynamical}.

The goal of this paper is to generate realizations of the stochastic process $s$ given its whitening operator $\mathrm{L}$ and a statistical characterization of its innovation process $w$.  The computer generation of these signals can be of great interest to practitioners who wish to evaluate their reconstruction algorithms. We are thinking of works such as \cite{Amini2013Interpolator,amini_bayesian_2013,kamilov_mmse_2013,godsill2006bayesian}, where optimal estimators for interpolating and  denoising such processes have been derived.

A possible approach to generate realizations of $s$ would be to notice that, if  $\mathrm{L}$ is a differential operator such as $\mathrm{D} = \frac{\mathrm{d}}{\mathrm{d}t}$ or a polynomial in $\mathrm{D}$, then \eqref{eq::found} defines a stochastic differential equation (SDE) \cite{oksendal_stochastic_1989}. This becomes more apparent when notating $w$ with the alternative notation  $\mathrm{d}Z_t$, where $(Z_t)_{t \in \mathbb{R}^+}$ is a Lévy process (Chapter 7.4 in \cite{unser_introduction_2014}). For example, $(\mathrm{D} - \alpha\mathrm{I})s = w$ can be rewritten as
$\mathrm{d}S_t = \alpha S_t \mathrm{d}t + \mathrm{d}Z_t$. A suitable SDE solver, such as the one studied in \cite{rubenthaler_numerical_2003}, can then be used to generate an approximation of the signal.  In particular, a common method is to solve the linear system of stochastic difference equations that is obtained by  considering the discrete counter-part of the operator ${\rm L}$ ({\it e.g.} using finite differences instead of the derivative), and by replacing  the innovation process $w$ with  a discrete white noise (see, for example, \cite{fageot_gaussian_2018}).

It turns out that generic SDE solvers do not exploit the linearity of $\mathrm{L}$. Here, the analytic treatment of \eqref{eq::found} can be pushed further to obtain an explicit solution. Brockwell  shows in \cite{brockwell_Levy-driven_2001} that $s$  corresponds to the integral of a deterministic function with respect to a Lévy process.  The integral can then be approximated by substituting it with a Riemann sum defined on a partition of the integration interval   \cite[Theorem 21]{protter_stochastic_2005}. 

These approaches, although valid, have drawbacks when it comes to the generation of synthetic signals for the evaluation of algorithms. First, they directly depend on the existence of a  grid on which the approximation of the continuous process is sampled. This can lead to complication in the context of the multi-resolution algorithms that  manipulate grid-free descriptions of signals. Second, the generated approximations are not solutions of an SDE in the form of \eqref{eq::found}. In other words, the approximations are not mathematical objects of the same nature as $s$.

In what follows, we propose a method that addresses both  issues. It is based on a theoretical result by  Fageot {\it et al.} \cite{fageot_gaussian_2018} that states that any solution $s$ of  \eqref{eq::found} is the limit in law of a sequence of simpler processes $s_n$. In other words, we have that 
\begin{equation*}
s_n \xrightarrow[]{\mathcal{L}} s, \text{ as } n \xrightarrow{} \infty.
\end{equation*}
These simpler processes, called generalized Poisson processes \cite{unser_introduction_2014}, have the advantage of having a grid-free numerical representation despite having a continuously defined domain. They fall within the category of (random) signals with a finite rate of innovation  \cite{vetterli_sampling_2002,unser2010stochastic}. They also have the desirable property of being whitened by the same operator $\mathrm{L}$ as the approximated signal. This implies that they all have the same correlation structure as the target signal (see Proposition \ref{prop1}). 

Our method takes a  sufficiently large  value for $n$ and generates a realization of the process $s_n$ on a chosen interval. To do so, we consider an intermediary process called the generalized increment process. Interestingly, this process can be represented  as a weighted sum of shifted B-splines and  can be sampled very efficiently \cite{unser_cardinal_2005, unser_cardinal_2005-1}. The desired stochastic process $s_n$ is then obtained from the latter by recursive filtering.

The outline of the paper is as follows: In Section II, we provide the necessary mathematical background. In Section III, we give a description of our algorithm: we begin by discussing the simulation of the innovation process in Subsection \ref{innov}. We then define the generalized increment process in Subsection \ref{subsec:genInc} and we show how to generate its trajectories in Subsection \ref{increment}. Using this, we provide a recipe for generating     sparse stochastic processes in Subsection \ref{subsec:recipe}. In Subsection \ref{subsec:corr}, we show that our generation method perfectly reproduces the correlation structure of the target stochastic process. Finally, we conduct numerical investigations to show the validity of our method in Section IV.

\section{Mathematical Foundations}

In this section, we give a brief overview of the mathematical concepts that underly  our approach. For a more detailed exposition, the reader is referred to \cite{unser_unified_2014,unser_unified_2014_2,fageot_gaussian_2018}, and references therein.

The Schwartz space $\Swa$ is the space of smooth and rapidly decaying test functions. Its continuous dual, denoted by $\Sp$, is the space of tempered distributions. It is the space of all continuous linear functionals over $\Swa$. 

We denote by $\mathrm{L}$ an  operator  that is a continuous, linear, shift-invariant mapping  from $\Sp$ to $\Sp$.  The operator $\rm L$ is said to be shift-invariant  if  for any test function $\varphi$ and any $t_0 \in\mathbb{R}$, we have that 
$$ {\rm L}\{ \varphi\} (t-t_0) = {\rm L}\{\varphi(\cdot-t_0) \}(t), \quad t\in\mathbb{R},$$
where $\varphi(\cdot -t_0): t \mapsto \varphi(t-t_0)$ is the shifted version of $\varphi$ by $t_0$.

 We restrict ourselves  to rational operators in $\mathrm{D} = \frac{\mathrm{d}}{\mathrm{dt}}$, written $\mathrm{L} = {P(\mathrm{D})}{Q(\mathrm{D})}^{-1}$, where $P$ and $Q$ are polynomials such that $\mathrm{deg}(P) > \mathrm{deg}(Q)$. The latter assumption is crucial to have the minimum required regularity (point-wise definition)  for the   solution $s$ of \eqref{eq::found}. The case ${\rm L}={\rm D}$ is a typical choice that   appears, for example, in the modeling of Brownian motion. 

 Rational operators are defined through their frequency response
\begin{equation*}
\widehat{\mathrm{L}}(\omega) = \frac{P(\mathrm{j}\omega)}{Q(\mathrm{j}\omega)}.
\end{equation*}
They provide a succinct representation of the equation $P(\mathrm{D})s = Q(\mathrm{D})w$ that we can simply rewrite as $\Lop s = w$.

We are interested in generalized stochastic processes defined over $\Sp$. A generalized stochastic process $w$ can be viewed as a random element of $\Sp$  in the sense that, for any $\varphi \in \Swa$, the linear functional $\varphi \mapsto \langle \varphi, w \rangle \in \mathbb{R}$ is a well defined random variable over $\mathbb{R}$  (See Appendix A for a formal definition).

\subsection{Lévy White Noises}\label{subsec:whitenoise}
 Lévy white noises  constitute an important class of generalized stochastic processes, whose specification is essential to our framework. The three important operational properties of L\'evy white noises for our purpose are: 
\begin{enumerate}
    \item  \textbf{Stationarity:}  For any $\varphi \in \Swa$ and $\tau \in \mathbb{R}$, the random variables $\langle \varphi, w \rangle$ and $\langle \varphi(\cdot - \tau), w \rangle$ are identically distributed.
    \item  \textbf{Independence:}  For any $\varphi_1, \varphi_2 \in \Swa$ with disjoint supports, the random variables $\langle \varphi_1, w \rangle$ and $\langle \varphi_2, w \rangle$ are independent.
    \item \textbf{Characterization of the probability law:} For any Lévy white noises $w$ in $\Sp$ and for any test function $\varphi \in \Swa$, the characteristic function of the random variable $X_{\varphi} = \langle \varphi, w\rangle$ can be specified as
\begin{equation}
\widehat{\mathcal{P}}_{X_{\varphi}}(\xi) =  \E[\mathrm{e}^{\mathrm{j}\xi\langle \varphi, w \rangle}] =\exp\left(\int_{\mathbb{R}} f(\xi \varphi(r))\mathrm{d}r\right),
\label{eq:charvar}
\end{equation}
where the function $f : \mathbb{R} \rightarrow \mathbb{C}$ is called the {Lévy exponent} of $w$.
\end{enumerate}

 Formally, this {Lévy exponent} can be obtained as 
$$f(\xi) = \log \left( \widehat{\mathcal{P}}_{X_{\textrm{rect}}}(\xi) \right),$$ 
where $X_{\textrm{rect}} = \langle \textrm{rect}_{[0, 1]}, w \rangle $\footnote{ Although   $\textrm{rect}_{[0, 1]}$ is not in $\Swa$,   the random variable $X_{\rm rect}$ can still be defined. For more details, see \cite{fageot_unified_2017}.} is the observation of $w$ through the rectangular window 
\begin{equation*}
\textrm{rect}_{[0, 1]}(x) = \begin{cases}
1, \quad  0 < x \leq 1\\
0, \quad \text{otherwise}.
\end{cases}
\end{equation*}
  The distribution of ${X_{\textrm{rect}}}$ gives us the {Lévy exponent} $f$ that defines \eqref{eq:charvar}, so that we can determine all the statistics of $w$  from the knowledge of  ${X_{\textrm{rect}}}$. 
  
 In particular, the following Proposition from \cite{unser_introduction_2014} connects the second-order statistics  of $w$ to those of ${X_{\textrm{rect}}}$. 
\begin{prop}[\cite{unser_introduction_2014}, Theorem 4.15]\label{prop1}
Let $w$ be a Lévy white noise such that $X_{\mathrm{rect}} = \langle \mathrm{rect}_{[0, 1]}, w \rangle$ has zero mean and a finite variance $\sigma_w^2 = \E[X_{\mathrm{rect}}^2]$. Then,
\begin{equation*}
\forall \varphi_1, \varphi_2 \in \Swa, \quad \E[\langle \varphi_1, w\rangle \langle \varphi_2, w\rangle ] = \sigma_w^2 \langle \varphi_1, \varphi_2 \rangle.
\end{equation*}
\end{prop}

It turns out that $ {X_{\textrm{rect}}}$  is an {infinitely divisible} random variable in the      sense of Definition \ref{Def:InfDiv} \cite{amini_sparsity_2014}. 

\begin{mydef}\label{Def:InfDiv}
A real-valued random variable $X$ is said to be infinitely divisible if, for any natural number $M \in \mathbb{N}$, there exist $M$ independent and identically distributed random variables $X_1, ..., X_M$ such that
\begin{equation*}
X = X_1 + \cdots + X_M.
\end{equation*}
\end{mydef}
To check the infinite divisibility of ${X_{\textrm{rect}}}$, one can note that, for any $M \in \mathbb{N}$, we have that
\begin{align}
{X_{\textrm{rect}}} = \langle \textrm{rect}_{[0, 1]}, w \rangle & = \langle \sum_{m = 0}^{M-1} \textrm{rect}_{[\frac{m}{M}, \frac{m+1}{M}]}, w \rangle \nonumber \\
& = \sum_{m = 0}^{M-1}\langle \textrm{rect}_{[\frac{m}{M}, \frac{m+1}{M}]}, w \rangle \label{eq:line2}.
\end{align}
The terms in the sum \eqref{eq:line2} are independent and identically distributed random variables as a consequence of the {independence} and {stationarity} properties of white noises, which certifies that $\langle \textrm{rect}_{[0, 1]}, w \rangle$ is infinitely divisible. 

The converse is also true: for any regular\footnote{ The random variable $X$ is said to be regular, if $\mathbb{E}[|X|^\epsilon]<+\infty$ for some $\epsilon>0$.} infinitely divisible random variable $X$ with Lévy exponent $f(\xi) = \log{(\E[\mathrm{e}^{\mathrm{j}\xi X}])}$, there exists a well defined Lévy white noise $w$ whose statistics are determined by \eqref{eq:charvar} \cite{amini_sparsity_2014,fageot_continuity_2014,dalang2017levy}. This shows that there is a one-to-one correspondence between infinitely divisible distributions and \textit{Lévy} white noises through $\langle \textrm{rect}_{[0, 1]}, w \rangle$.

The Gaussian, gamma, and $\alpha$-stable distributions are classical examples of infinitely divisible distributions \cite{sato}. We can plug in their Lévy exponents in \eqref{eq:charvar} to define their corresponding Lévy white noises. We repeat in Table \ref{tab:inftab}  some infinitely divisible distributions of interest,  along with their Lévy exponents\cite{amini_sparsity_2014}.
\renewcommand{\arraystretch}{1.43}
\begin{table}[t]
\centering
\caption{Infinitely divisible distributions and their Lévy exponents}
\begin{tabular}{l|l}
\hline\hline
 \textbf{Distribution}     &  \textbf{Lévy exponent}      \\  \hline
Gaussian $(\mu, \sigma)$             & $\mathrm{j} \mu \xi  - \sigma^2 \xi^2/2$      \\ 
Symmetric $\alpha$-stable $(\alpha,  c), \alpha \in (0,2]$ & $ -|c\xi|^{\alpha}$                      \\  
Gamma$(\alpha, \beta)$               & $-\beta \log{(1 - {\mathrm{j}\xi}/{\alpha})}$ \\  
Laplace  $(\mu, b)$      & $\mathrm{j}\mu \xi  - \log{(1  + b^2\xi^2)}$         \\  \hline \hline
\end{tabular}
\label{tab:inftab}
\end{table}

A case of special interest  is when $f$ is the Lévy exponent of a compound-Poisson distribution. A compound-Poisson random variable $X$, with rate $\lambda$ and amplitude law $\nu$, is defined as
\begin{equation*}
X = \sum_{k=1}^{K}A_k,
\end{equation*}
where the number $K$ is a Poisson random variable with parameter $\lambda$ and $(A_k)_{k =1}^K$ is an i.i.d. sequence drawn according to $\nu$. We refer to the corresponding Lévy white noise $w$ as a  compound-Poisson innovation. It is known to be equal in law to
\begin{equation}\label{Eq:CPNoise}
w = \sum_{k \in \mathbb{Z}} A_k \delta( \cdot - \tau_k),
\end{equation}
where $(\tau_k)_{k \in \mathbb{Z}}$ are the locations of impulses with rate $\lambda$ \cite{unser_unified_2014}. The law of these impulses is as follows: for any interval $[a, b]$, the number of impulses in $[a,b]$ is a Poisson random variable with parameter $\lambda(b-a)$. 

On any finite interval, compound-Poisson innovations have a finite representation. They can be stored on a computer with the quantization of real numbers as sole source of information loss. They are therefore well adapted to simulation purposes.

\subsection{Generalized Lévy Processes}\label{subsec:genLevy}

The sparse-stochastic-process framework of Unser {\it et al.} \cite{unser_introduction_2014} is a comprehensive theory of {generalized Lévy Processes}. These are stochastic processes that can be whitened by some admissible linear, shift-invariant operator. More precisely, $s$ is a generalized Lévy process if there exists an operator $\Lop$ such that $w = \Lop s$
is a Lévy white noise. Equivalently, one may view generalized Lévy processes as the solution of the stochastic differential equation
\begin{equation}
    \Lop s = w.
\label{eq:lsw}
\end{equation}
It has been shown that, under mild technical assumptions on $\mathrm{L}$ and $w$, a solution $s$ of \eqref{eq:lsw} exists and constitutes a properly defined generalized stochastic process over $\Sp$ \cite{fageot_continuity_2014}. 

When $\Lop$ is an operator with a trivial null space, such as $\Lop=(\mathrm{D} - \alpha \mathrm{I})$ with $\Re(\alpha) \neq 0$, we can write that
$$
s = \Lop^{-1}w,
$$
where $\Lop^{-1}$ is the inverse of $\Lop$. However, when the null space is nontrivial, for instance when $\Lop$ corresponds to an unstable ordinary differential equation, the specification of the boundary conditions become necessary to uniquely identify the solution. The boundary conditions  take the form
\begin{equation}
\phi_\ell(s)= c_\ell, \quad \ell=1,\ldots,N_0,
\label{boundary}
\end{equation}
where $\phi_\ell: s\mapsto \phi(s)\in\mathbb{R}$ are appropriate linear functionals, $c_\ell \in \mathbb{R}$, and $N_0$ is the dimension of the null space of $\mathrm{L}$. For instance, one can impose that  the process $s$ takes fixed values at reference locations $t_1<\ldots<t_{N_0}$; that is, $\phi_\ell(s)=s(t_\ell)=c_{\ell}$ for $\ell=1,\ldots,N_0$.  Such boundary conditions appear in the classical definition of Lévy processes (including Brownian motion), where we have that $\phi(s) = s(0) =  0$ (Chapter 7 of  \cite{unser_introduction_2014}). We formally write 
\begin{equation*}
s = \mathrm{L}_{\boldsymbol{{\phi}}}^{-1}w,
\end{equation*}
where $\mathrm{L}_{\boldsymbol{{\phi}}}^{-1}$ is the right inverse of $\mathrm{L}$.  It incorporates the boundary conditions \eqref{boundary} (Chapter 5.4 of \cite{unser_introduction_2014}).

When $w$ is a compound-Poisson innovation of the form \eqref{Eq:CPNoise}, the process $s = \mathrm{L}_{\boldsymbol{{\phi}}}^{-1}w$ ($\mathrm{L}^{-1}s$, respectively,  when the null space of $\mathrm{L}$ is trivial) is called a {generalized Poisson process}.  

The fundamental property for this work is that any sparse stochastic process $s$ that is the solution of \eqref{eq:lsw} can be specified as the limit in law of a sequence $\{ s_n \}_{n \in \mathbb{N}}$  of generalized Poisson processes   \cite{fageot_gaussian_2018}. The corresponding driving processes $w_n = \mathrm{L}s_n$ are compound-Poisson innovations of the form
\begin{equation}\label{Eq:WnForm}
w_n = \sum_{k \in \mathbb{Z}} A_{k,n} \delta( \cdot - \tau_{k,n})
\end{equation}
with rates $\lambda_n = n$ and with i.i.d.\ amplitudes $A_{k,n}$ that are infinitely divisible random variables with Lévy exponent $f_n = \frac{1}{n}f$, where $f$ is the Lévy exponent of $w$.

\subsection{Green's Functions}\label{subsec:greens}

The Green's  function of a differential operator $\Lop$ is a tempered distribution $\rho_\Lop \in \Sp$ that satisfies
\begin{equation*}
\Lop\rho_\Lop= \delta.
\end{equation*}
It can be viewed as the impulse response of the inverse of $\Lop$. The canonical Green's function is
\begin{equation*}
\rho_\Lop =  \mathcal{F}^{-1}\left\{ \frac{1}{\widehat{\Lop}(\omega)}\right\},
\end{equation*}
where $\widehat{\Lop}$ is the frequency response of $\Lop$ (Chapter 5.2 of \cite{unser_introduction_2014}). This definition can be made to stay valid even when $\widehat{\Lop}$ vanishes at some points, as long as $\frac{1}{\widehat{\Lop}(\omega)}$ is in $\Sp$.  For details on how to compute Green's functions, the reader is referred to   Appendix B. We have plotted the Green's function of several operators  in Figure \ref{fig:BsplineGreens}  to highlight their variety and their dependence on ${\rm L}$.

% Note that Green's functions are not unique. They can differ by any element of the null-space of $\Lop$. For instance for the first order differential operator $\Lop = \mathrm{D} $, $x \mapsto \frac{1}{2}\mathrm{sign}(x) - c$ is a Green's function for any $c\in\mathbb{R}$. However, only $\rho_{\mathrm{D}, 0}(x) = \frac{1}{2}\mathrm{sign}(x) - \frac{1}{2}$ satisfies the boundary condition 
% $$
% \langle \varphi, \rho_{\mathrm{D}, 0} \rangle := \int_{-\infty}^{0}\rho_{\mathrm{D}, 0}(x)\mathrm{d}x = 0.
% $$

% This example illustrates that in general there exists a unique Green's function $\rho_L$ of $\Lop$ that satisfies boundary conditions of the form
% \begin{equation}
% \langle \varphi_\ell, \rho_{\Lop} \rangle := \int_{-\infty}^{0}\rho^{(\ell-1)}_{\Lop, 0}(x)\mathrm{d}x = 0, \text{ $\ell \in \{1, ..., N_0 \}$},
% \label{eq::boundary}
% \end{equation}
% where $N_0$ is the dimension of the null space of $\Lop$. 

\section{Method}

In this section, we introduce our method for generating (approximate) trajectories of a sparse stochastic process $s$ that is whitened by an operator $\Lop$ and whose innovation noise is $w$. When necessary, we assume general boundary conditions of the form $\phi_\ell(s)= 0$ for $\ell = 1, ..., N_0$, where $N_0$ is the dimension of the null space of  $\Lop$.

As mentioned earlier, the process $s$ is the limit of generalized compound-Poisson processes $s_n$ driven by $w_n = \Lop s_n$, a compound-Poisson innovation of the form \eqref{Eq:WnForm}. The process $s_n$ can therefore be written 
\begin{equation*}
s_n =  \sum_{k \in \mathbb{Z}} A_{k,n} \rho_{\Lop}( \cdot - \tau_{k,n}) + p_{0,n},
\end{equation*}
where $\rho_{\Lop}$ is a Green's function of $\Lop$ and $p_{0,n}$  is an element of the null space of $\Lop$ determined by boundary conditions (it vanishes when $\Lop$ is invertible). Indeed, we have that
\begin{align*}
\Lop \left\{ \sum_{k \in \mathbb{Z}} A_{k,n} \rho_{\Lop}( \cdot - \tau_{k,n}) + p_{0,n}\right\} & = \sum_{k \in \mathbb{Z}} A_{k,n} \Lop \{\rho_{\Lop}( \cdot - \tau_{k,n})\} \nonumber\\
& = \sum_{k \in \mathbb{Z}} A_{k,n} \delta( \cdot - \tau_{k,n}) \nonumber \\&= w_n.
\end{align*}
For large values of $n$, the process $s_n$ is assumed to be a good approximation of $s$. So, our goal is to generate samples of $s_n$ on any uniform grid over any interval $[0, T]$.  More precisely, once an interval $[0, T]$ is specified and a regular grid with step size $h$ is provided, our aim is to obtain the vector $\mathbf{s}_n$ whose components are $[\mathbf{s}_n]_i = s_n(ih)$, for $i = 0, ..., \left(\lceil \frac{T}{h} \rceil - 1\right)$.
\subsection{Simulating the Innovation Process}
\label{innov}
We begin by obtaining a realization of the driving innovation $w_n$. It consists of a sequence of impulse locations $(\tau_{k,n})$ and a corresponding sequence of amplitudes $(A_{k,n})$.

The sequence $(\tau_{k,n})$ is a point Poisson process. Its realization on the interval $[0, T]$ is simulated in two steps. First, a Poisson random variable $K$ with parameter $\lambda = nT$ is generated. Then, $K$ impulse locations $(\tau_{k,n})_{k \in \{1,\ldots,K\}}$ are sampled uniformly on $[0, T]$.

The next step is to simulate the $K$ corresponding amplitudes $(A_{k,n})_{k \in \{1,\ldots,K\}}$. The characteristic function of the amplitudes variable $A$ is
\begin{equation*}
\xi \mapsto \exp \left(\frac{1}{n}f(\xi)\right).
\end{equation*}
We refer to it as the \textit{n}th root of the law of $\langle \textrm{rect}_{[0, 1]}, w \rangle$. Our assumption 
in this paper is that there exists, for any $n \in \mathbb{N}$, a known method\footnote{Workarounds exists for when a sampling method for 
the $n$th root $\frac{1}{n}f(\xi)$ is unavailable. For instance, one can opt for an approximate sampling scheme such as in 
\cite{bondesson_simulation_1982}.} to generate infinitely divisible 
variables with Lévy exponent $\frac{1}{n}f(\xi)$.    For common parametric distributions such as $\alpha$-stable, Laplace, and gamma 
distributions, such sampling methods\cite{devroye_chapter_2006} are well known and implemented in scientific 
computing libraries\footnote{E. Jones,  et al., “SciPy: Open source scientific tools for Python,” 2001.}. Simulating from their \textit{n}th root is a simple matter 
of rescaling their parameters, as summarized in Table \ref{tab:my-table}. By applying the correct rescaling, we simulate $K$ independent amplitudes and thus obtain the sequence $(A_{k,n})_{k \in \{1,\ldots,K\}}$. 
\renewcommand{\arraystretch}{1.7}
\begin{table}[t]
\centering
\caption{The $n$th root of infinitely divisible distributions}
\begin{tabular}{l|l}
\hline\hline
\textbf{Distribution}  & \textbf{$n$th Root}                                                                                         \\ \hline
Gaussian $(\mu, \sigma)$                                & Gaussian $(\frac{\mu}{n}, \frac{\sigma}{\sqrt{n}})$                                                               \\ 
$\alpha$-Stable $(\alpha, \beta, \mu, c)$  & If $\alpha \neq 1$, $(\alpha, \beta, \frac{\mu}{n}, \frac{c}{n^{\frac{1}{\alpha}}})$,                           \\
                                                        & If $\alpha = 1$, $(\alpha, \beta, \frac{\mu}{n} - \frac{2}{\pi}c\beta\frac{\log(n)}{n}, \frac{c}{n})$ \\ 
Gamma$(\alpha, \beta)$                                  & Gamma$(\frac{\alpha}{n}, \beta)$                                                                                        \\ 
Compound-Poisson of intensity $\lambda$                        & Compound-Poisson of intensity $\frac{\lambda}{n}$                                                                             \\ 
Laplace $(\mu, b)$                                      & $X_n = \frac{\mu}{n} + b(G^{(n)}_1 - G^{(n)}_2)$ \\
 & with $G^{(n)}_1, G^{(n)}_1 \sim \textrm{Gamma}(\frac{1}{n}, 1)$ \\
 \hline\hline

\end{tabular}
\label{tab:my-table}
\end{table}

\subsection{ Generalized Increment Process}\label{subsec:genInc}

With the impulse locations $(\tau_{k,n})_{k \in \{1,\ldots,K\}}$ and amplitudes $(A_{k,n})_{k \in \{1,\ldots,K\}}$ in hand, we    can compute samples of 
\begin{equation}
\label{eq:greenrepr2}
s_{n}(\cdot)= \sum_{k =1}^{K} A_{k,n} \rho_{\Lop}( \cdot - \tau_{k,n}) + p_{0,n}
\end{equation}
on a grid. 

 A direct approach to generate $\mathbf{s}_n$ is to use the expansion \eqref{eq:greenrepr2} and  represent the process as a sum of shifted Green's functions. However  in this case, the determination of $s_n(t)$ at any point $t \in [0, T]$ may require nontrivial computation  of each and every term in \eqref{eq:greenrepr2}. This stems from the fact that Green's functions are infinitely supported in general. There are therefore potential drawbacks to expansions in the basis of shifted Green's functions like \eqref{eq:greenrepr2}. To overcome these issues, we propose instead an alternative method based on B-splines.

\begin{figure}[t]
\centering
\includegraphics[ width=0.5\textwidth]{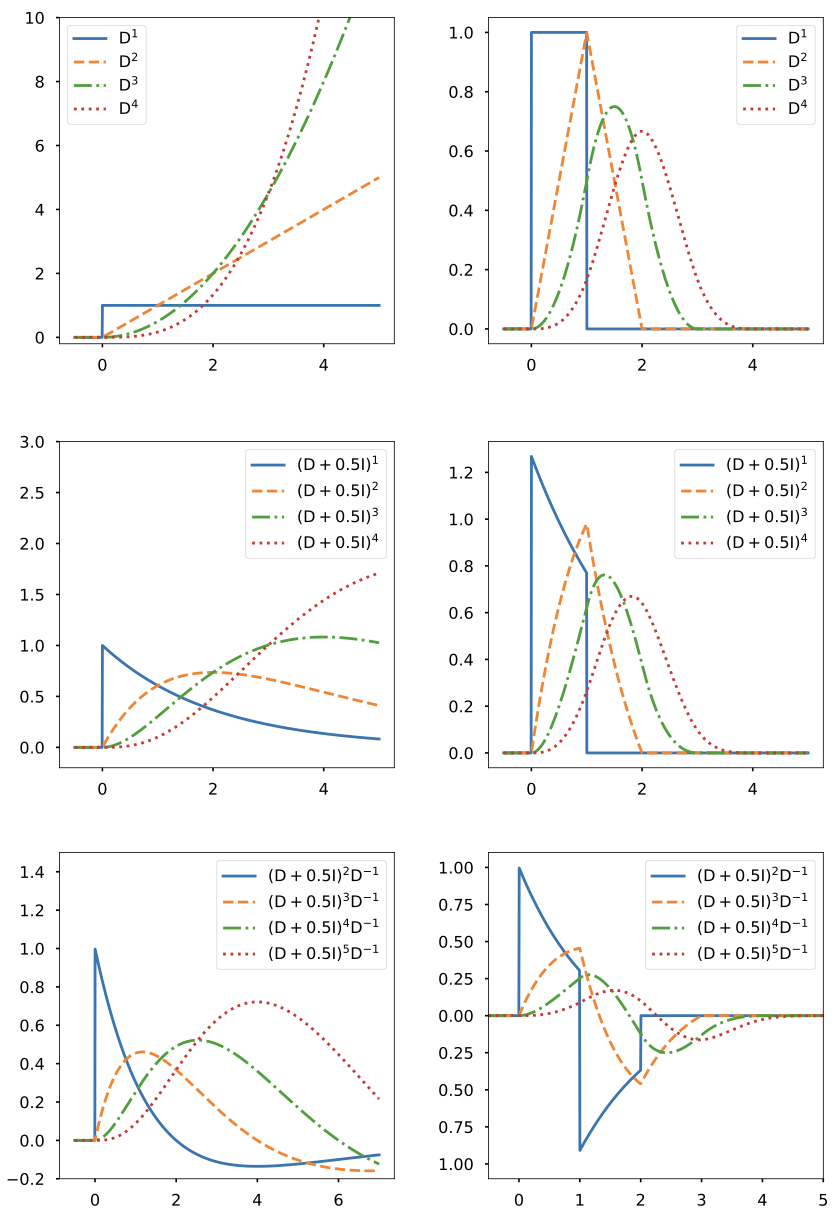}
\caption{Green's functions (left) and B-splines (right) associated with   specific operators $\rm L$.}
\label{fig:BsplineGreens}
\end{figure}

 Recall that $\Lop$ is a rational operator of the form ${P(\mathrm{D})}{Q(\mathrm{D})}^{-1}$, where we take $\{\alpha_1, ..., \alpha_{\mathrm{deg}(P)} \}$ to be the roots of $P$, with possible repetitions. Its discrete counterpart $\Lop_{\rm d}^h$ is defined as
\begin{equation*}
\Lop_{\rm d}^h\{f\} = \sum_{m=0}^{\mathrm{deg}(P)}r[m] f(\cdot - mh),
\end{equation*}
where the sequence $r$ is determined through its Fourier transform
$$
 R({\rm e}^{{\rm j}\omega}) = \sum_{m=0}^{{\rm deg}(P)} r[m] {\rm e}^{-{\rm j}\omega m} = \prod_{m=1}^{\mathrm{deg}(P)}(1 - {\rm e}^{\alpha_m h}{\rm e}^{-{\rm j}\omega h}).
$$
It is a finite impulse-response filter (FIR). Its null space contains the null space of $\Lop$ \cite{unser_cardinal_2005}. The function  $\beta_{\rm L}^h := \Lop_{\rm d}^h\{\rho_{\Lop}\}$ is called the \mbox{B-spline} corresponding to $\Lop$ \cite{de1972calculating}. The B-spline has the fundamental property of being the shortest possible function within the space of cardinal ${\rm L}$-splines (its support is included in $[0, \textrm{deg}(P)\times h]$) \cite{schoenberg1988contributions,ron1990factorization}. This will turn out to be crucial for the numerical efficiency of our method. Moreover, they reproduce both the Green's function  and  elements in the null space of their corresponding operator ${\rm L}$ \cite[Section 6.4.]{unser_introduction_2014}. Examples of relevant  generalized B-splines are shown in Figure \ref{fig:BsplineGreens} (right figures). Note how they contrast with the corresponding infinitely supported Green's functions (left figures).
 
The application of $\Lop_{\rm d}^h$ to $s_n$ yields  
\begin{align} 
u_n(t) = \Lop_{\rm d}^h\{s_{n}\}(t) = \sum_{m=0}^{\mathrm{deg}(P)}r[m] s_n(t - mh).
\label{eq:udef}
\end{align}
The process $u_n$ in \eqref{eq:udef} is called the generalized increment process. Interestingly, it can be written as a sum of compactly supported terms, like 
\begin{align*} 
u_n(t)& = \Lop_{\rm d}^h\{s_{n}\}(t)\\ & = \sum_{k=1}^{K} A_{k,n}
 \Lop_{\rm d}^h\{\rho_{\Lop}( \cdot - \tau_{k,n})\}(t)  + \Lop_{\rm d}^h\{p_{0,n}\}(t) \nonumber\\
& = \sum_{k=1}^{K} A_{k,n} \beta_{\rm L}^h (t - \tau_{k,n}) + \quad 0. 
\end{align*}
The process $u_n$, along with boundary conditions, is our alternate representation of $s_n$. Now, let $\mathbf{u}_n$ be the vector whose components are $[\mathbf{u}_n]_i = u(ih)$, for $i = 1, ...,  \left(\lceil\frac{T}{h} \rceil - 1\right)$. This vector can be computed more efficiently than $\mathbf{s}_n$ since the process $u_n$ admits a representation with compactly supported terms. Moreover,   ${\bf u}_n$ is linearly related to the vector ${\bf s}_n$ via a discrete system of difference equations. Indeed, we have that 
\begin{equation}
[\mathbf{u}_n]_i = \sum_{m=0}^{\mathrm{deg}(P)}r[m] [\mathbf{s}_n]_{i - m}, 
\label{eq::filter}
\end{equation}
for $\mathrm{deg}(P) \leq i \leq \left(\lceil \frac{T}{h} \rceil - 1\right)$. For $ 0 < i < \mathrm{deg}(P)$, we have that 
\begin{equation*}
[\mathbf{u}_n]_i = \sum_{m=0}^{\mathrm{deg}(P)}r[m] \mathbf{s}_n((i - m)h),
\end{equation*}

where the values  $\mathbf{s}_n(-mh)$ for $m = 0, ..., \left(\mathrm{deg}(P)-1\right)$  provide the boundary values. These relations are established by writing \eqref{eq:udef} with $t = ih$. The boundary values are determined by the null-space term $p_{0,n}$, which is itself determined by the boundary conditions. 

Thus, once we have evaluated $\mathbf{u}_n$,  we can obtain $\mathbf{s}_n$ by solving \eqref{eq::filter}, which is accomplished by applying a recursive reverse filter to $\mathbf{u}_n$. This  is performed by rewriting  \eqref{eq::filter} as  
\begin{equation}
[\mathbf{s}_n]_i = \frac{1}{r[0]} \left( [\mathbf{u}_n]_i - \sum_{m=1}^{\mathrm{deg}(P)}r[m] [\mathbf{s}_n]_{i - m} \right).
\label{eq::reverse}
\end{equation}
By substitution of the boundary values when necessary ( {\it i.e.},  taking $\mathbf{s}_n((i - m)h)$ instead of $[\mathbf{s}_n]_{i - m}$ when $(i-m) \leq 0$), \eqref{eq::reverse} allows one to recursively compute the components of $\mathbf{s}_n$.

\subsection{Computing the Generalized Increment Process}
\label{increment}
We now describe an efficient procedure to compute the generalized increment process. The components of $\mathbf{u}_n$ are given by
\begin{equation*}
[\mathbf{u}_n]_i = \sum_{k=0}^{K} A_{k,n} \beta_{\rm L}^h (ih - \tau_{k,n}).
\end{equation*}
The naive approach here would be to iterate through each grid point $i$ independently and compute [$\mathbf{u}_n]_i$. Doing so would require one to read the entire sequence of impulse locations $(\tau_k)$ for each $i$. This cannot be avoided since there is no information on the sequence $(\tau_k)$, aside from its inclusion in $[0,T]$. We simply would not know \textit{which} B-spline terms are inactive, so we would have to iterate through them all. A more efficient approach is to  iterate through the list of impulses instead of  the grid points. 
 \begin{figure}[t]
    \centering
    \includegraphics[ width=0.4\textwidth]{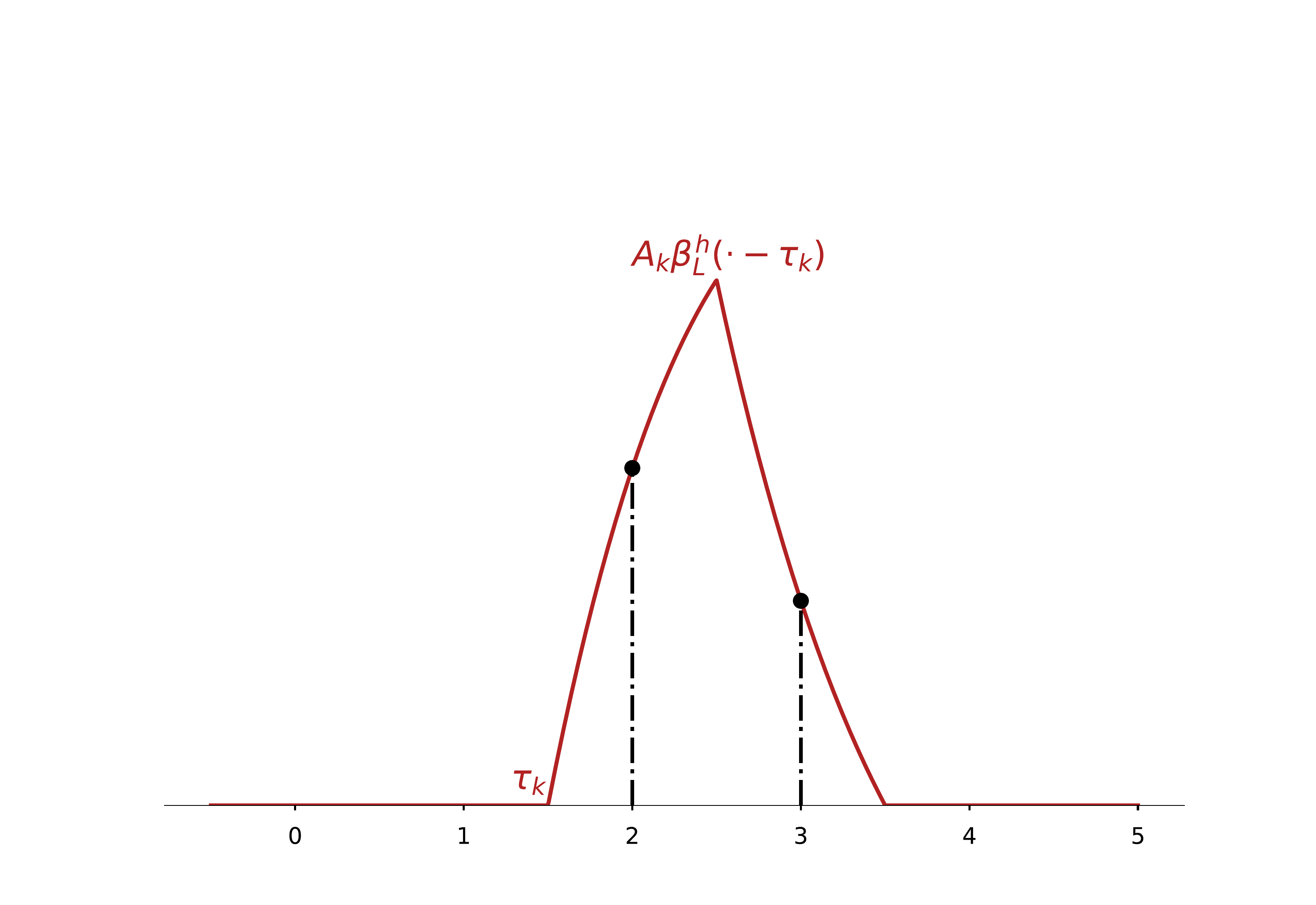}
    \caption{For a single B-spline term, it is only the grid points that sit within the support of the B-spline that are incremented (black stems).}
    \label{fig::illustration}
\end{figure}

The idea is as follows: First, initialize the vector $\mathbf{u}_n$ to zeros. Then, read the list of impulse locations one by one. For each impulse at $\tau_k$, find the grid points that lie within the support of the B-spline at $\tau_k$. Then, increment the value of $\mathbf{u}_n$ on those grid points by the contribution of the considered B-spline (see Figure \ref{fig::illustration}). In one pass over the list of impulses, this method computes the values $[\mathbf{u}_n]_i = \mathbf{u}_n(ih)$.

This intermediate computation of the generalized-increment process provides a considerable gain in terms of efficiency. Instead of having a number of operations that scales with $\lceil \frac{T}{h} \rceil \times K$ for the Green's function representation, we have one that scales with $\mathrm{deg}(P) \times \left(\lceil \frac{T}{h} \rceil + K\right)$.

\subsection{Recipe to  Generate Trajectories}\label{subsec:recipe}

Here is a summary of  the procedure that generates trajectories of $\Lop s_n = w_n$.

First, fix the infinitely divisible distribution\footnote{The choice here is restricted to parametric families we can rescale and simulate.} that corresponds to $w$  and define the operator $\Lop$ by identifying the polynomials $P$ and $Q$.

Pick a sufficiently large value for $n$. Intuitively, $n$ should be large enough to ensure the occurrence  of several jumps in each bin. In other words, we expect $n$ to be of the same order as $h^{-1}$. This has been validated with our numerical experiments as well, where we show that it provides a good approximation of the underlying statistics of the process  (see Subsection \ref{Subsec:ConvN} and Figures \ref{fig:browninan}, \ref{fig:convergencetwo}, and \ref{fig:KS}). 

Pick a simulation interval $[0, T]$ and generate $w_n$ as described in section \ref{innov}. Determine an explicit form for $\rho_\Lop$. At this point, the grid-free approximation $s_n$ (expressed as in \eqref{eq:greenrepr2}) is available and can be stored. 

Fix a grid on $[0, T]$ by choosing a step size $h$. Then  determine  the vector $\mathbf{s}_n$ with component  $[\mathbf{s}_n]_i= s_n(ih)$. Compute the FIR filter $\Lop_{\rm d}^h$ and obtain $\beta_\Lop^h = \Lop_{\rm d}^h \{\rho_\Lop\}$. Then, compute the generalized increment vector $\mathbf{u}_n$ as described in Section \ref{increment}. 

To obtain $\mathbf{s}_n$, apply the reverse filter to $\mathbf{u}_n$ following \eqref{eq::reverse}. Take  the values $\mathbf{s}_n(-mh)$ for $m = 0, ..., \left(\mathrm{deg}(P)-1\right)$ to be zero for most cases except when $\Lop$ has a nontrivial null space, in which case it is derived from boundary conditions. The pseudocode of our method is provided in Algorithm \ref{algo}.

\begin{algorithm}[t]
 
 \SetKwInOut{Input}{Input}
\SetKwInOut{Output}{Output}
 \Input{Coefficients of $P$ and $Q$, approximation level $n$, interval size $T$, step size $h$}
 \Output{Vector $\mathbf{s}_n$}
 Compute $\rho_\Lop$ and the FIR filter $r[m]$ \\ 
 Compute $\beta_L^h= \Lop_d^h\{\rho_{\Lop}\}$ \\
 Generate [$(\tau_{1}, A_{1}), ..., (\tau_{K}, A_{K})$] \\
 Initialize  $\mathbf{u}_n$ with zeros as an array of size $\lceil \frac{T}{h} \rceil$ 
 
 \ForEach{$(\tau_k, A_k)$}
 {
   Find closest grid point $i_{{\rm grid}} = \lfloor \frac{\tau_i}{h}\rfloor$ \\
   \ForEach{i in $\{i_{{\rm grid}} , \ldots, i_{{\rm grid}} + \mathrm{deg}(P) \}$}{
   $[\mathbf{u}_n]_i \xleftarrow[]{} [\mathbf{u}_n]_i +  A_k \times \beta_L^h(ih - \tau_k)$  
   }
 }
Recursively apply a reverse filter to $\mathbf{u}_n$ following \eqref{eq::reverse} 
\vspace{2mm}
  \caption{Procedure to obtain $\mathbf{s}_n$.}\label{algo}
\end{algorithm}

\subsection{Correlation Structure}\label{subsec:corr}
In this section, we  show a merit of our method by proving that the generated approximations preserve the correlation structure of the target process.

First  note that for any white Lévy noise $w$, we have that 
\begin{equation*}
w_n \xrightarrow[]{\mathcal{L}} w,
\end{equation*}
where the sequence of compound-Poisson innovations $(w_n)_{n \in \mathbb{N}}$ is defined  in \eqref{Eq:WnForm}. We refer to this approximating sequence in Proposition \ref{Prop}.

\begin{prop}\label{Prop}
Let $w$ be a Lévy white noise such that $X_{\mathrm{rect}} = \langle \mathrm{rect}_{[0, 1]}, w \rangle$ has zero mean and the finite variance $\sigma_w^2 = \E[X_{\mathrm{rect}}^2]$. Let $n \in \mathbb{N}$  and let $w_n$ be a compound-Poisson innovation that approximates  $w$ as defined in \eqref{Eq:WnForm}. Denoting $X_{\mathrm{rect}, n} =  \langle \mathrm{rect}_{[0, 1]}, w_n \rangle$, we have that 
\begin{equation*}
\E[X_{\mathrm{rect}, n}] = \E[X_{\mathrm{rect}}] = 0
\end{equation*}
and
\begin{equation*}
\sigma_{w_n}^2  = \E[X_{\mathrm{rect}, n}^2] = \E[X_{\mathrm{rect}}^2] = \sigma_w^2.
\end{equation*}
\end{prop}
The  proof can be found in Appendix C. Now, if $s_n = \Lop^{-1}w_n$ is a generalized Poisson process that approximates $s = \Lop^{-1}w$, then
\begin{align}
\E[\langle \varphi_1, s_n\rangle \langle \varphi_2, s_n\rangle ] & = \E[\langle \varphi_1, \Lop^{-1}w_n\rangle \langle \varphi_2, \Lop^{-1}w_n\rangle ] \nonumber \\
& = \E[\langle \Lop^{-1*}\varphi_1, w_n\rangle \langle \Lop^{-1*}\varphi_2, w_n\rangle ] \nonumber\\
& =  \sigma_{w_n}^2 \langle \Lop^{-1*}\varphi_1, \Lop^{-1*}\varphi_2 \rangle \nonumber \\
& =  \sigma_{w}^2 \langle \Lop^{-1*}\varphi_1, \Lop^{-1*}\varphi_2 \rangle \nonumber\\
& = \E[\langle \varphi_1, s\rangle \langle \varphi_2, s\rangle ]. \label{Eq:CorrStruct}
\end{align}
From \eqref{Eq:CorrStruct}, we concluded that, more than just approximated, the correlation structure is  preserved \textit{exactly} in our method. 
\section{Numerical Experiments}
In this section, we validate our approach by conducting several numerical experiments.  Let us also mention that a Python library that implements our algorithm can be found online\footnote{ \href{https://github.com/Biomedical-Imaging-Group/Generating-Sparse-Processes}{https://github.com/Biomedical-Imaging-Group/Generating-Sparse-Processes}}. Moreover, an accompanying web interface is also designed and is available \footnote{\href{https://saturdaygenfo.pythonanywhere.com}{https://saturdaygenfo.pythonanywhere.com}}.

\subsection{Generating Lévy Processes} 
\begin{figure}[t]
    \centering
    \includegraphics[width=0.5\textwidth]{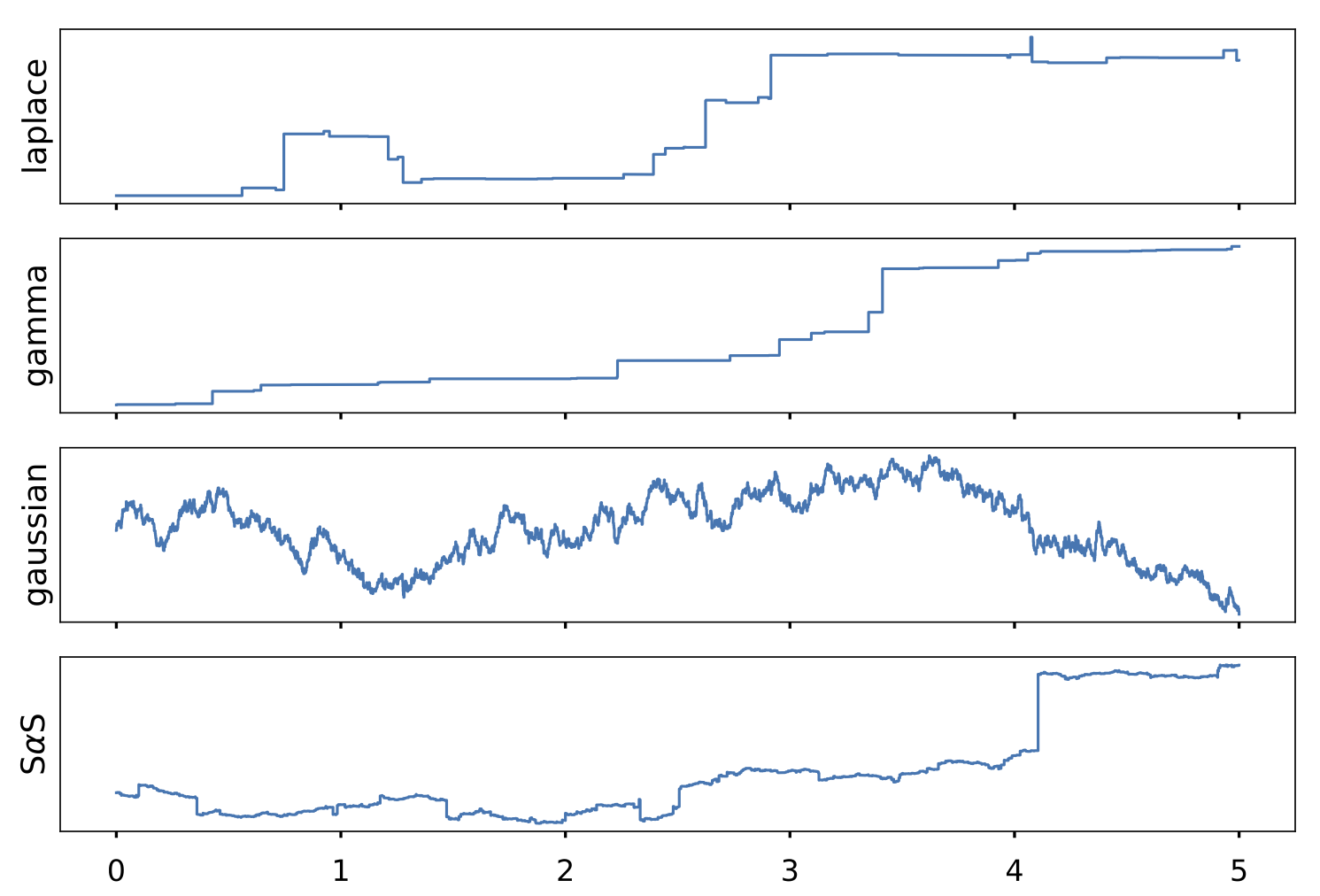}
    \caption{Trajectories of Lévy processes ($\Lop = \mathrm{D}$) with different innovations. From top to bottom: Laplace(0, 1), gamma(1, 1),  Gaussian(0,1), and symmetric-$\alpha$-stable with $\alpha = 1.23$.}
    \label{fig:diffw}
\end{figure}

Among all processes we can generate, those that are solutions to $\mathrm{D}s = w$ are called Lévy processes when the boundary condition is $s(0) = 0$. We showcase in Figure \ref{fig:diffw} different Lévy processes that correspond to several infinitely divisible distributions. For all four simulations, we took $n = 1,\!000$ and  $h = 0.001$. As we demonstrate in Section \ref{Subsec:ConvN}, a reasonable choice for these parameter is to set $nh$ to be a small integer (here,   $nh=1$).  The visual appearance of the trajectories matches our expectations: The trajectory driven by a Gaussian innovation has the appearance of Brownian motion; the gamma Lévy process is nondecreasing.

\subsection{Choice of  the Operator} 
\begin{figure}[t]
\centering
\includegraphics[width=0.5\textwidth]{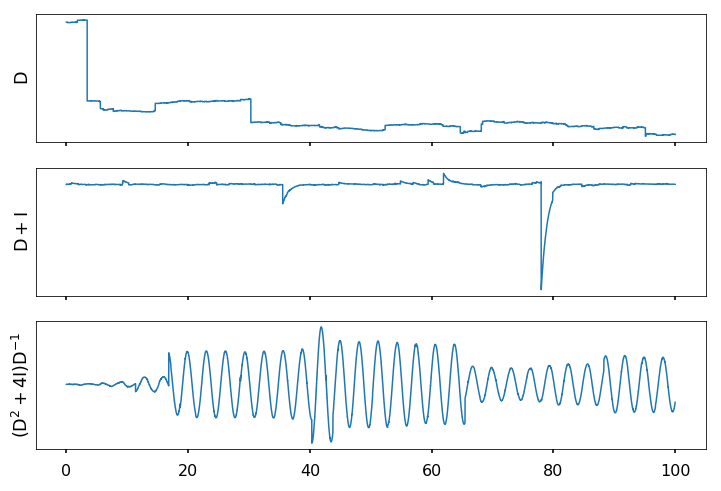}
\caption{ Trajectories of the solution $s$ of ${\rm L}s=w$ for different  operators $\rm L$. In all cases, we considered  a symmetric-$\alpha$-stable white noise $w$ with $\alpha=1.23$.}
\label{fig:diffop}
\end{figure}

Our framework allows for any rational operator of the form ${P(\mathrm{D})}{Q(\mathrm{D})}^{-1}$, so long as $\mathrm{deg}(P) > \mathrm{deg}(Q)$. In Figure \ref{fig:diffop}, we generate trajectories of $s$ that are solution of $\Lop s = w$, where $w$ is a symmetric-$\alpha$-stable innovation with $\alpha = 1.23$. Here we took $n = 200$ and $h = 0.001$. We  see that, for various choices of $\Lop$, the characteristics of the signal are markedly different, which exhibits the breadth of the modeling framework proposed in \cite{unser_introduction_2014}.

\subsection{ Convergence as $n$ Grows} 
\label{Subsec:ConvN}

 \begin{figure}[t]
    \centering
    \includegraphics[ width=0.5\textwidth]{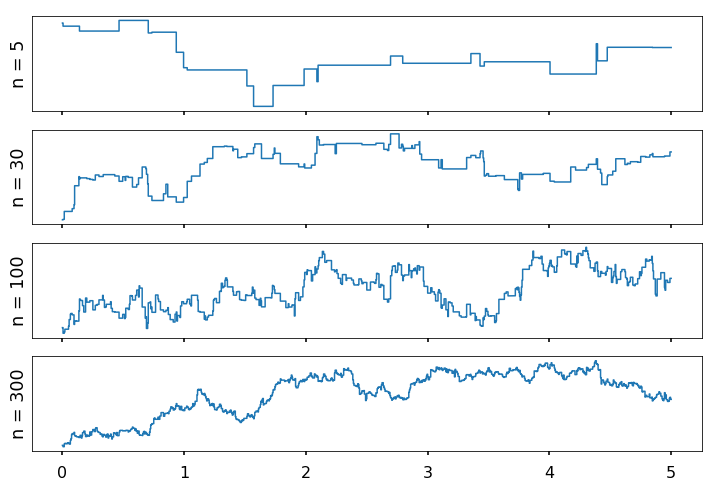}
    \caption{Approximations of Brownian motion (solution to $\Dop s = w$, with $w$ a Gaussian white Lévy noise) as $n$ increases.}
    \label{fig:browninan}
\end{figure}

 \begin{figure}[t]
    \centering
    \includegraphics[ width=0.5\textwidth]{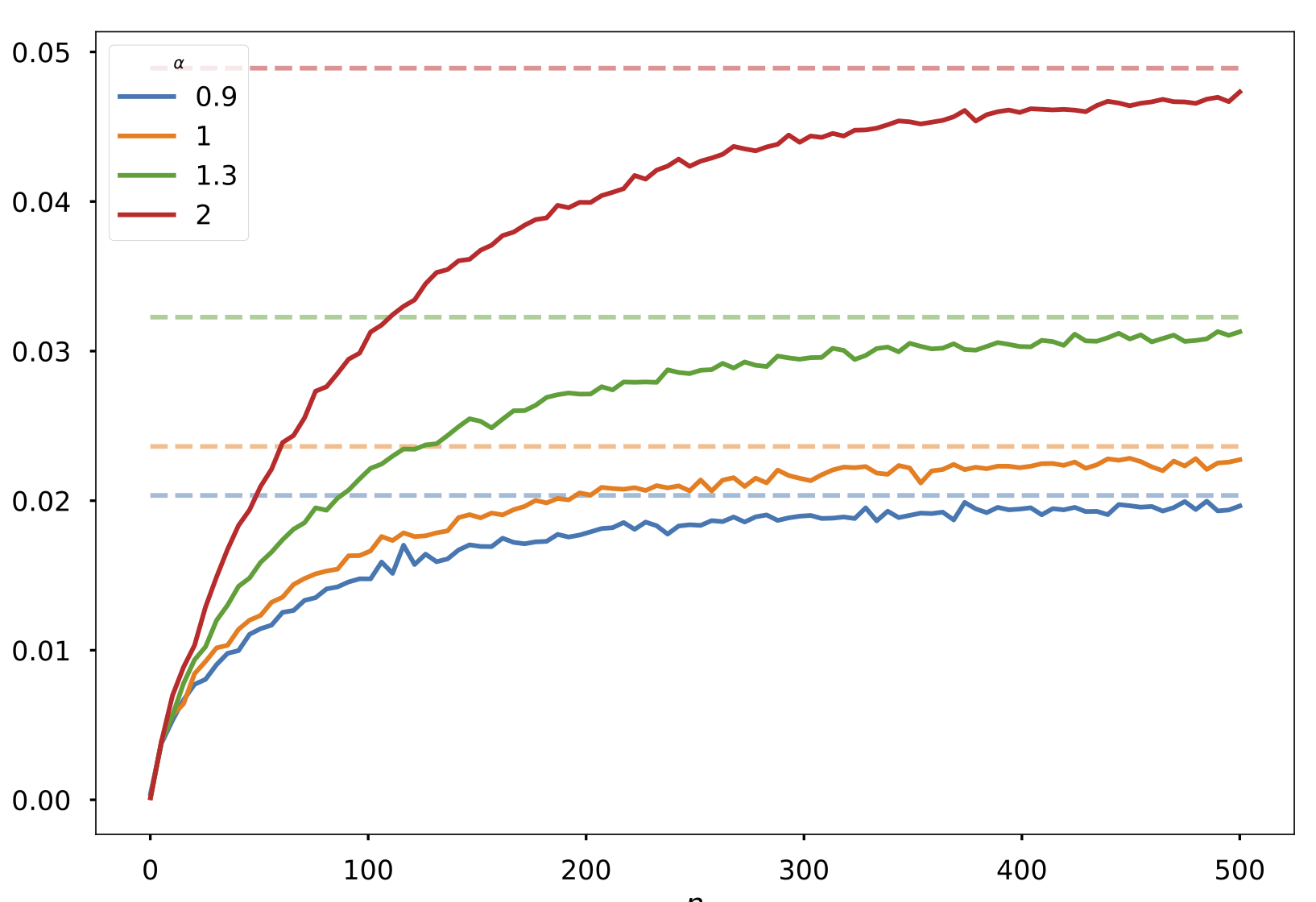}
    \caption{Convergence of $\E [ |\langle \mathrm{rect}_{[0, h]}, s_n \rangle|^{p}]$ to $\E [ |\langle \mathrm{rect}_{[0, h]}, s \rangle|^{p}]$ for $p=0.4$, $h=0.01$,  and several symmetric-$\alpha$-stable Lévy white noises $w$. The expectations are estimated with $10,\!000$ trajectories for each $n$.}
    \label{fig:convergencetwo}
\end{figure}

 \begin{figure}[t]
    \centering
    \includegraphics[width=0.5\textwidth]{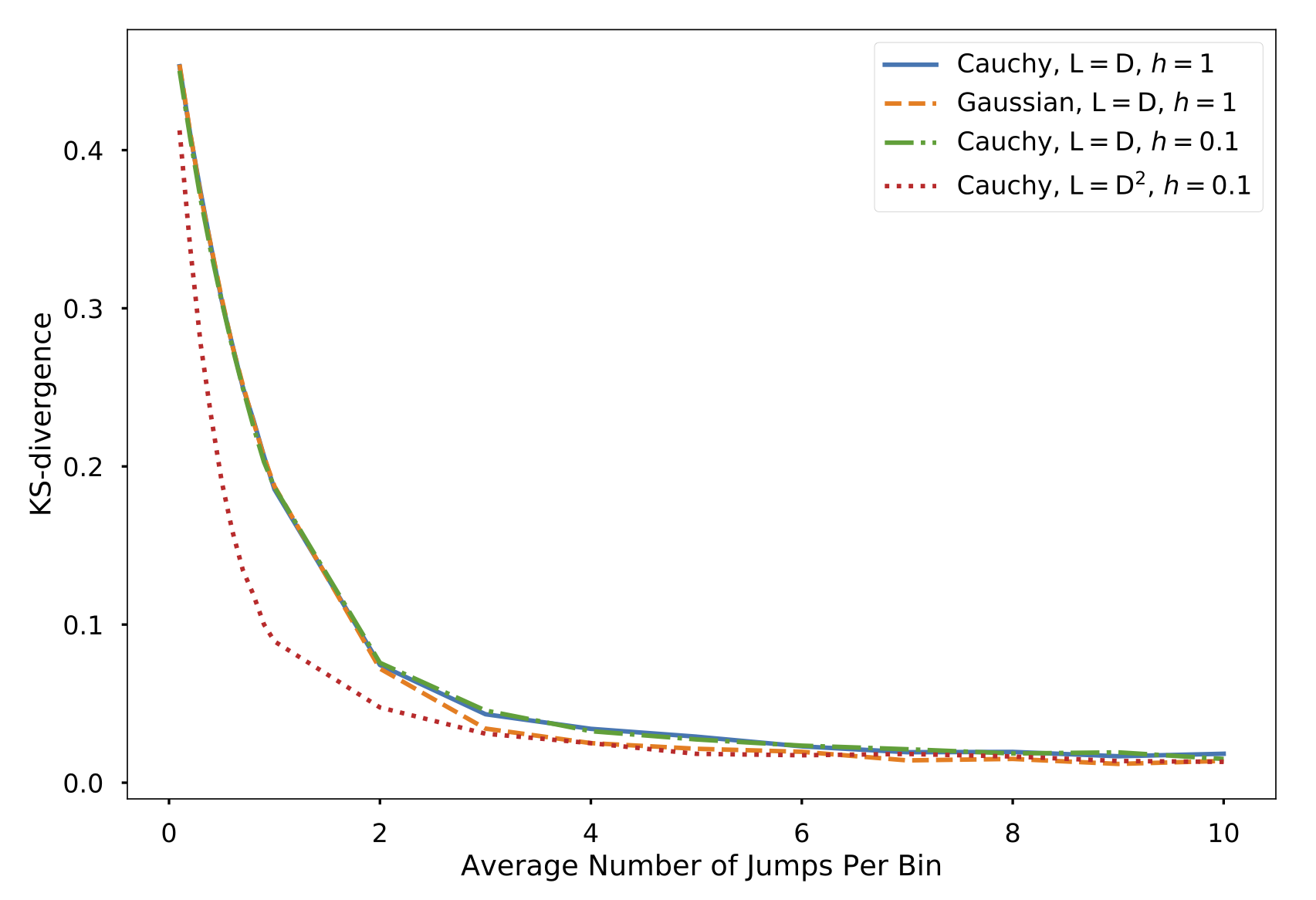}
    \caption{ Kolmogorov-Smirnov (KS) divergence versus the average number of jumps per bin ($N_{\text{jumps}}=nh$).   }
    \label{fig:KS}
\end{figure}
 In Figure \ref{fig:browninan}, we illustrate  how an increase in $n$ improves the approximation. In addition, we have depicted the convergence  of moments in Figure \ref{fig:convergencetwo}. While the two figures emphasize the effect of  $n$,  they are insufficient to provide a quantitative way to choose  $n$. 

Here, we propose a measure that is based on the statistics of the generalized increment process. Since the process $u_n$ is maximally decoupled, we can estimate the   distribution of $U_n= \langle \beta_{\mathrm{L}}^{h\vee}, w_n\rangle$  from the samples $\{[\mathbf{u}_n]_i\}_i$ of the generalized increment process on the grid and obtain the empirical cumulative  distribution function (CDF) $\tilde{F}_{n}(\cdot)$ of $U_n$. We then compare this empirical function to the reference CDF $F(\cdot)$ of $U= \langle \beta_{\mathrm{L}}^{h\vee}, w \rangle$. For the comparison, we use the Kolmogorov-Smirnov (KS) divergence \cite{massey1951kolmogorov}  defined as 
\begin{equation*}
\mathcal{KS}(\tilde{F}_n,F) = \max_{x\in \mathbb{R}} |\tilde{F}_n(x) - F(x) |.
\end{equation*}
We   then select $n$ such that the KS-divergence is smaller than a certain threshold ({\it e.g.}, smaller than $0.1$). The choice of the threshold is conditioned by the desired numerical precision: The lower the threshold, the more faithful the trajectories, but  the higher the computational cost of the algorithm. 

Intuitively, we expect that it is necessary to have several jumps in each bin  in order to { properly} approximate the statistics of the process. The average number of jumps in each bin  of length $h$  is $N_{\text{jumps}}=nh$, so we expect $n$ to be in the  order of $h^{-1}$.  

In Figure \ref{fig:KS}, we have validated this intuition by plotting the \mbox{KS-divergence} for different values of $N_{\text{jumps}}$ in various settings. In all cases, as $N_{\text{jumps}}$ increases, the KS-divergence decreases to a baseline error value,   due to the finite-sample estimation of the underlying distribution.

\subsection{Benefits of Grid-Free Approximations} 
\begin{figure}[t]
    \centering
    \includegraphics[ width=0.5\textwidth]{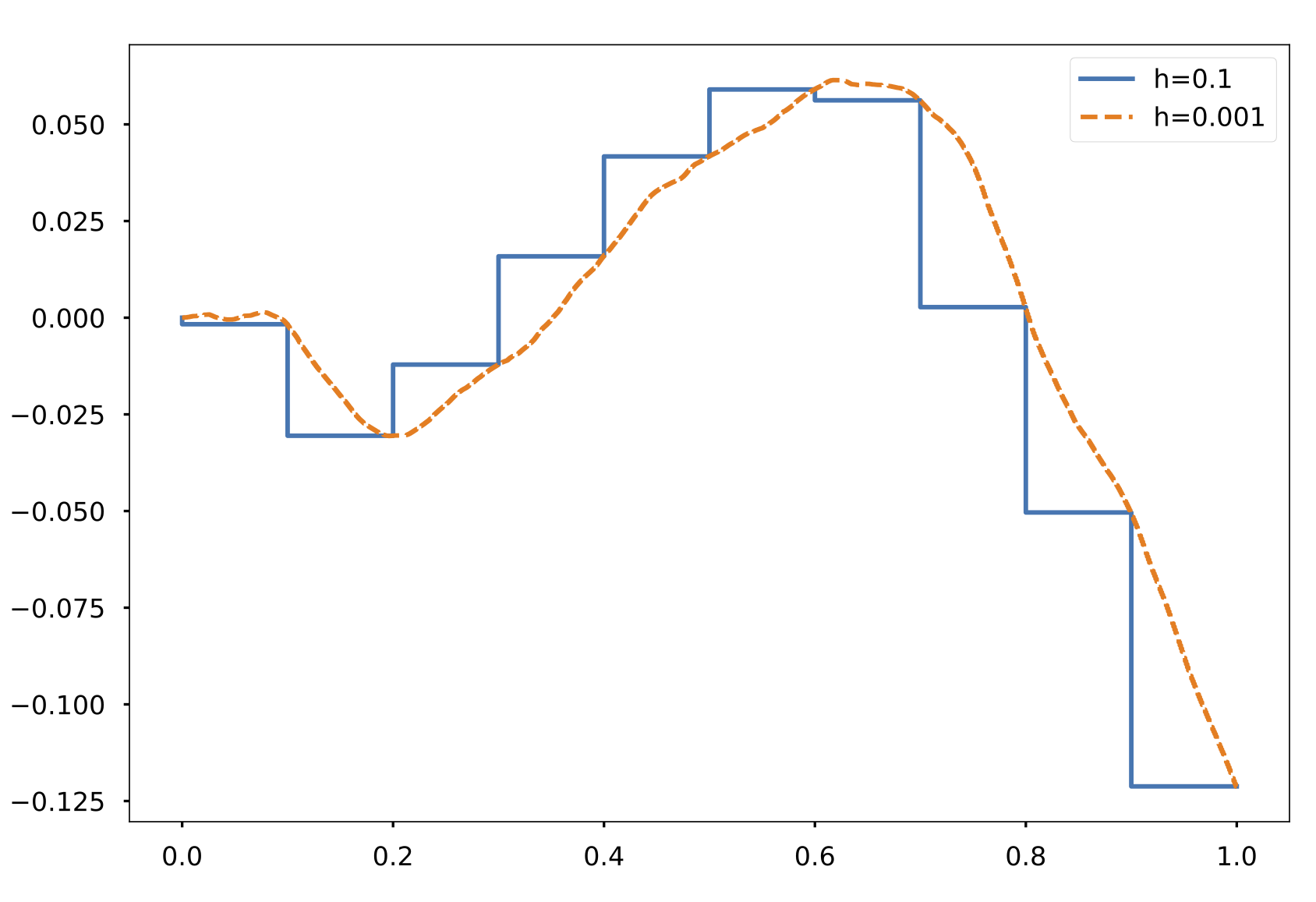}
    \caption{Single grid-free approximation sampled on  grids that differ by their step size.}
    \label{fig:multigrid}
\end{figure}

Recall that a main motivation for our algorithm was to make it   compatible with multi-grid methods. In our approach, the approximation $s_n$ lives off the grid. It is only after the specification of the step size $h$ that $s_n$ is sampled on a grid. The generation of the random variables to determine $s_n$ and the sampling on a grid are completely decoupled. This means that the same  approximation $s_n$ can be viewed through different grids, which we illustrate  in Figure \ref{fig:multigrid}. The solution to $(\Dop+1)^2s = w$, where $w$ is a Gaussian white Lévy noise, is first approximated by $s_{1000}$. Then, it is viewed on  different regular grids on $[0, 1]$.

\subsection{Computational Efficiency} 
\begin{figure}[t]
    \centering
    \includegraphics[ width=0.5\textwidth]{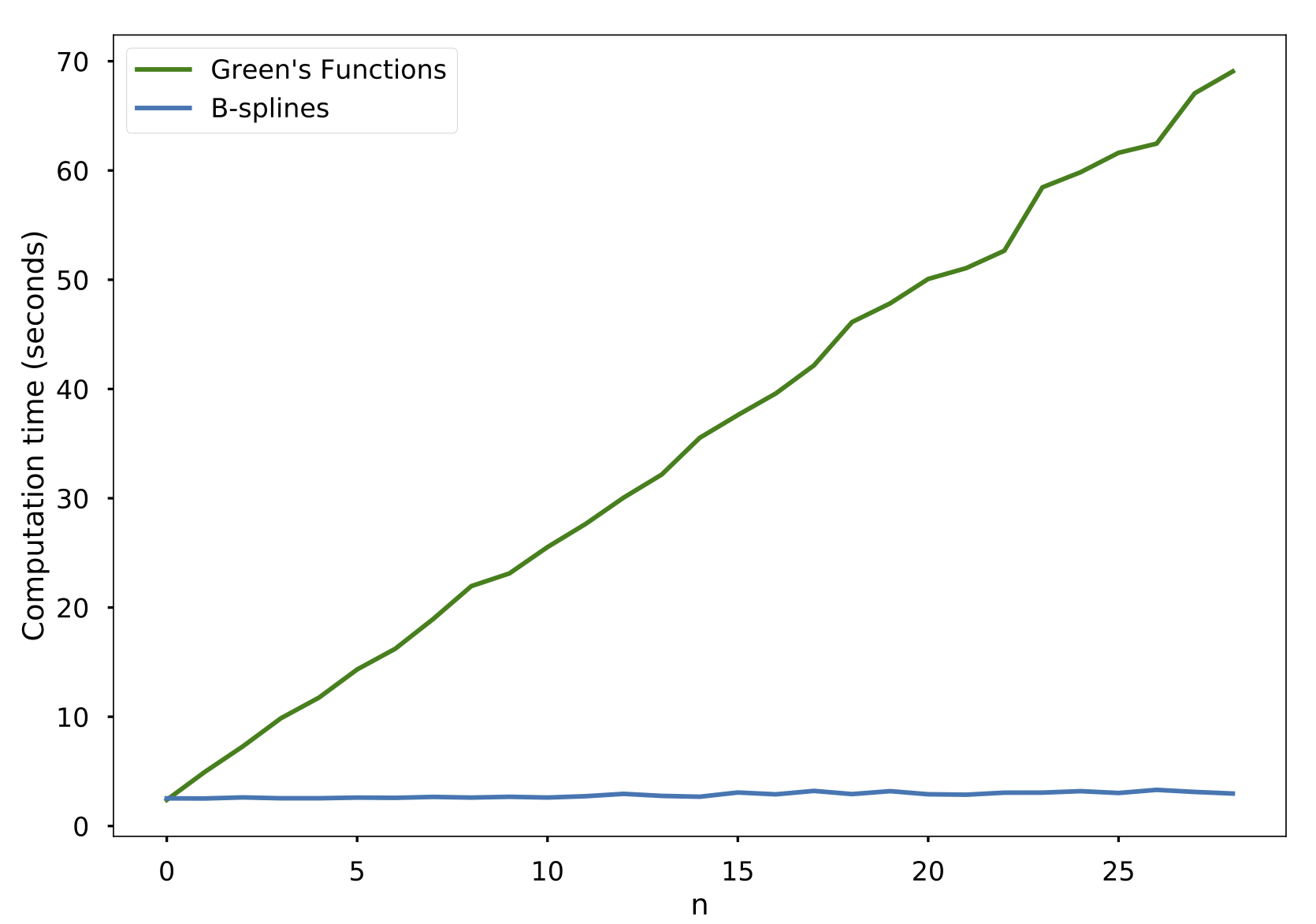}
    \caption{Average computation time for a trajectory of the solution of $(\mathrm{D} - 0.5\mathrm{I})s = w$, where $w$ a Gaussian white noise. The simulation interval is $[0, 1]$ with step size $h = 0.001$.}
    \label{fig:computime}
\end{figure}

A crucial component of our approach is the computation of the generalized increment   $\mathbf{u}$ in order to obtain the values of $s_n$ on a grid. This provides a gain in numerical efficiency that can be felt even on moderately sized simulations. As can be seen in Figure \ref{fig:computime}, using a Green's function representation requires significantly more time than using an intermediate B-spline representation.

\section{Conclusion}

We have described a novel approach for generating sparse stochastic processes. Our method leverages the properties of B-splines to guarantee good numerical efficiency.  A possible direction for future work is to provide theoretical guidance on how one should choose  the  parameter $n$ in terms of a prescribed  tolerance on the approximation error.

\section*{Acknowledgment}

The authors would like to thank Dr. Julien Fageot, Pakshal Bohra, and Thomas Debarre for enlightening discussions.

% if have a single appendix:
%\appendix[Proof of the Zonklar Equations]
% or
%\appendix  % for no appendix heading
% do not use \section anymore after \appendix, only \section*
% is possibly needed

% use appendices with more than one appendix
% then use \section to start each appendix
% you must declare a \section before using any
% \subsection or using \label (\appendices by itself
% starts a section numbered zero.)
%
\appendices

\section{Generalized stochastic processes}

Generalized stochastic processes are random elements of $\Sp$ that can be fully specified by their characteristic functionals. Those are infinite-dimensional generalizations of the characteristic functions of real random variables.
\begin{mydef}
The characteristic functional of the generalized stochastic process $s$ is the functional $\widehat{\mathcal{P}}_s : \Swa \xrightarrow{} \mathbb{C}$ such that
\begin{equation*}
\widehat{\mathcal{P}}_s(\varphi) = \E[{\rm e}^{{\rm j}\langle \varphi, s \rangle}], \text{ for all }\varphi \in \Swa.
\end{equation*}
It is a continuous, positive-definite functional and $\widehat{\mathcal{P}}_s(0) = 1$.
\end{mydef}
Just as in finite dimensions, $\widehat{\mathcal{P}}_s$ contains all the statistical information of $s$. In particular, for any test function $\varphi \in \Swa$, the distribution of the real random variable $\langle \varphi, s \rangle$ is entirely determined by $\widehat{\mathcal{P}}_s$ as its probability density function $p$ is given by
\begin{equation*}
p(t) \propto  \mathcal{F}^{-1} \left\{\E[{\rm e}^{{\rm j}\omega\langle \varphi, s \rangle}]\right\}(t) =  \mathcal{F}^{-1} \left\{ \widehat{\mathcal{P}}_s(\omega\varphi)\right\}(t),
\label{eq:var}
\end{equation*}
where $\mathcal{F}^{-1}$ is the inverse Fourrier transform.
The construction of such objects was initiated in  \cite{gelfand_generalized_2014}. Their use for modeling sparse signals was introduced in \cite{unser_introduction_2014}.

\section{Computing Green's functions}

Here, we describe a method to compute Green's functions of rational operators. We begin with the intermediate computation of the Green's function of $\Lop = (\mathrm{D} - \alpha \mathrm{I})^k$. We  have that 
\begin{align}
\rho_{\alpha, k}(t) & = \mathcal{F}^{-1} \left \{ \frac{1}{ ({\rm j}\omega - \alpha )^k}\right \} (t) \nonumber\\
& = \begin{cases}
\mathbbm{1}_+(t) \frac{t^{k - 1}}{(k - 1)!}{\rm e}^{\alpha t} ,  & \Re(\alpha) \leq 0\\
-\mathbbm{1}_+(-t) \frac{t^{k - 1}}{(k - 1)!}{\rm e}^{\alpha t}, &  \text{ otherwise}
\end{cases}
\label{eq:basicelement}
\end{align}
is a Green's function of $\Lop$. 

Now, recall that rational operators are of the form $\mathrm{L} = {P(\mathrm{D})}{Q(\mathrm{D})}^{-1}$, where $P$ and $Q$ are polynomials. Taking $\{\alpha_1, ..., \alpha_{m} \}$ to be the roots of $P$ with multiplicity $\{\gamma_1, ..., \gamma_{m} \}$, the inverse of the frequency response is given by
\begin{equation*}
\frac{1}{\widehat{\Lop}(\omega)} = \frac{Q({\rm j}\omega)}{\prod_{i=1}^{m}({\rm j}\omega - \alpha_i)^{\gamma_i}} .
\end{equation*}
This inverse is known to admit a partial-fraction decomposition of the form
$$
\frac{1}{\widehat{\Lop}(\omega)} = \sum_{i=1}^{m}\sum_{k=1}^{\gamma_i}\frac{c_{ik}}{({\rm j}\omega - \alpha_i)^k}
$$
for some constants $c_{ik} \in \mathbb{C}$. The corresponding Green's function is then given by:
\begin{align*}
\rho_\Lop (t) & =  \mathcal{F}^{-1}\left \{ \frac{1}{\widehat{\Lop}(\omega)}\right \}  (t) \nonumber \\
& = \sum_{i=1}^{m}\sum_{k=1}^{\gamma_i}c_{ik}\mathcal{F}^{-1}\left \{ \frac{1}{({\rm j}\omega - \alpha_i)^k}\right \} (t) \nonumber \\
& = \sum_{i=1}^{m}\sum_{k=1}^{\gamma_i}c_{ik} \rho_{\alpha_i, k} (t)
\end{align*}
 The Green's function of $\Lop$ is then be obtained by summing   the Green's function of the partial fractions given in \eqref{eq:basicelement}.

\section{}
Proof of Proposition 1: Since $w_n$ is a compound-Poisson innovation,  $X_{\textrm{rect, n}}$ is a compound-Poisson random variable. It can be written
\begin{align*}
X_{\textrm{rect,n}} = \sum_{i=1}^{N}A_i
\end{align*}
where $N$ is a Poisson random variable with rate $\lambda = n$ and the $A_i$  are independent identically distributed infinitely divisible random variables with Lévy exponent $\frac{1}{n}f$ independent from $N$. We have by independence of the $(A_i)$, that 
\begin{equation*}
X_{\textrm{rect,n}}=\sum_{i=1}^{n}A_i =_d X_{\textrm{rect}}
\end{equation*}
because the characteristic function of $\sum_{i=1}^{n}A_i$ is $({\rm e}^{\frac{1}{n}})^n ={\rm e}^f$. This directly implies that $X_{\textrm{rect,n}}$ and $X_{\textrm{rect}}$ have the same moments. 
%\begin{align*}
%n\E[A_1]= \E [ \sum_{i=1}^{n}A_i ] & = \E[X_{\textrm{rect}}] =0.
%\end{align*}
%Moreover, 
%\begin{align*}
% \E[X_{\textrm{rect}}^2] & = \E [ (\sum_{i=1}^{n}A_i)^2 ] \\
% & =  \E [ \sum_{i=1}^{n}A_i^2 ] \quad \text{(by independence)} \\
% & = n \E[A_1^2]. 
%\end{align*}
%We can now show the equality of moments. We only show the equality of second order moments as the derivations are similar:
%\begin{align*}
%\E [ X_{\textrm{rect, n}}^2 ] & = \E [ (\sum_{i=1}^{N}A_i)^2]  =\E [ \E[ (\sum_{i=1}^{N}A_i )^2|N ]]\\
%& = \sum_{k=0}^{+\infty} \E[\sum_{i=1}^{k}A_i^2]] \P(N = k) \\
%& = \E[A_1^2]\sum_{k=0}^{+\infty} k \P(N = k) \\
%& =\E[A_1^2] \E [N ]  \\
%& = \frac{\E[X_{\textrm{rect}}^2]}{n} n = \E[X_{\textrm{rect}}^2]
%\end{align*}

% Can use something like this to put references on a page
% by themselves when using endfloat and the captionsoff option.
\ifCLASSOPTIONcaptionsoff
  \newpage
\fi

% trigger a \newpage just before the given reference
% number - used to balance the columns on the last page
% adjust value as needed - may need to be readjusted if
% the document is modified later
%\IEEEtriggeratref{8}
% The "triggered" command can be changed if desired:
%\IEEEtriggercmd{\enlargethispage{-5in}}

% references section

% can use a bibliography generated by BibTeX as a .bbl file
% BibTeX documentation can be easily obtained at:
% http://mirror.ctan.org/biblio/bibtex/contrib/doc/
% The IEEEtran BibTeX style support page is at:
% http://www.michaelshell.org/tex/ieeetran/bibtex/
%\bibliographystyle{IEEEtran}
% argument is your BibTeX string definitions and bibliography database(s)
%\bibliography{IEEEabrv,../bib/paper}
%
% <OR> manually copy in the resultant .bbl file
% set second argument of \begin to the number of references
% (used to reserve space for the reference number labels box)
\bibliographystyle{IEEEtran.bst}
\bibliography{lvy.bib}
% biography section
% 
% If you have an EPS/PDF photo (graphicx package needed) extra braces are
% needed around the contents of the optional argument to biography to prevent
% the LaTeX parser from getting confused when it sees the complicated
% \includegraphics command within an optional argument. (You could create
% your own custom macro containing the \includegraphics command to make things
% simpler here.)
%\begin{IEEEbiography}[{\includegraphics[width=1in,height=1.25in,clip,keepaspectratio]{mshell}}]{Michael Shell}
% or if you just want to reserve a space for a photo:

%\begin{IEEEbiography}{Michael Shell}
%Biography text here.
%\end{IEEEbiography}

% if you will not have a photo at all:
%\begin{IEEEbiographynophoto}{John Doe}
%Biography text here.
%\end{IEEEbiographynophoto}

% insert where needed to balance the two columns on the last page with
% biographies
%\newpage

%\begin{IEEEbiographynophoto}{Jane Doe}
%Biography text here.
%\end{IEEEbiographynophoto}

% You can push biographies down or up by placing
% a \vfill before or after them. The appropriate
% use of \vfill depends on what kind of text is
% on the last page and whether or not the columns
% are being equalized.

%\vfill

% Can be used to pull up biographies so that the bottom of the last one
% is flush with the other column.
%\enlargethispage{-5in}

% that's all folks
\end{document}